\documentclass[leqno]{amsproc}
\usepackage{amsfonts}
\usepackage{amssymb, latexsym, amsmath, pb-diagram}
\usepackage{pstricks-add}
\usepackage{multicol}
\usepackage{bm}
\usepackage[mathscr]{eucal}
\usepackage{xy}
\usepackage{epic,eepic}
\xyoption{all}

\def\w{\widetilde}

\def\o{\overline}

\newcommand{\ff}{\mathbb{F}}
\newcommand{\zz}{\mathbb{Z}}

\newcommand{\RR}{\mathbb{R}}

\begin{document}
\title[A Homeomorphism Invariant of Polyhedra]
{\bf A Homeomorphism Invariant of Polyhedra}
\author[Q. Zheng]{Qibing Zheng}
\thanks{The author is supported by NSFC grant No. 11071125}
\date{}
\address{School of Mathematical Science and LPMC, Nankai University,
Tianjin, 300071, P.R.China}
\email{zhengqb@nankai.edu.cn}
\maketitle

\begin{abstract} In this paper, we define a new bigraded $L$-homology
on finite simpicial complexes and prove that $L$-homology is a homeomorphism invariant of polyhedra.
\end{abstract}
\vspace*{4mm}

The following definition is a brief review of finite simplicial complex theory.\vspace{3mm}

{\bf Definition 1} A (finite, abstract) simplicial complex $K$ with vertex set $S$
is a set of subsets of the finite set $S$ satisfying the following two conditions.

1) $\phi\in K$ and for all $s\in S$, $\{s\}\in K$.

2) If $\sigma\in K$ and $\tau\subset\sigma$, then $\tau\in K$.

An element $\sigma$ of $K$ is called an $n$-simplex of $K$ if $\sigma$ has $n{+}1$ vertices of $S$ and $|\sigma|=n$ is the dimension of $\sigma$. Specifically, the empty simplex $\phi$ has dimension $-1$. A face of a simplex $\sigma$ is a subset $\tau\subset\sigma$. A proper face is a proper subset. A simplicial subcomplex $L$ of $K$ is a subset $L\subset K$ such that $L$ is also a simplicial complex. For a simplex $\sigma$ of $K$, the link of $\sigma$ is the simplicial subcomplex link$_{K}\sigma=\{\tau\in K\,|\,\sigma{\cup}\tau\in K,\,\sigma{\cap}\tau=\phi\}$.

For two simplicial complexes $K$ and $L$, their union is the simplicial complex $K*L=\{\sigma{\sqcup}\tau \,|\,\sigma{\in}K,\,\tau{\in}L\}$ ($\sqcup$ is the disjoint union). The cone $CK$ is the union of $K$ and the simplicial complex with only one vertex $v$.

For two simplicial complexes $K$ and $L$ with respectively vertex set $S$ and $T$, their Cartesian product is the simplicial complex $K\bar\times L$ with vertex set $S{\times}T$ defined as follows. Give $S$ and $T$ an order. For two non-empty simplexes $\sigma{=}\{v_1,\cdots,v_m\}\in K$ and $\tau{=}\{w_1,\cdots,w_n\}\in L$ such that $v_1{<}\cdots{<}v_m$ and $w_1{<}\cdots{<}w_n$, $\sigma\bar\times\tau$ is the set of subsets of $\sigma{\times}\tau$ of the form $\{(v_{i_1},w_{j_1}),\cdots,(v_{i_s},w_{j_s})\}$ such that $i_1{\leqslant}\cdots{\leqslant}i_s$ and $j_1{\leqslant}\cdots{\leqslant}j_s$. Then $K\bar\times L=\{\phi\}\cup_{\sigma\in K,\tau\in L,\sigma,\tau\neq\phi}\sigma\bar\times\tau$. Notice that by definition, if one of $K$ and $L$ is the empty set $\{\phi\}$, then $K\bar\times L=\{\phi\}$.

For simplicial complexes $K$ and $L$ with respectively vertex set $S$ and $T$, a simplicial map $f\colon K\to L$ is a map $f\colon S\to T$ such that for all $\sigma\in K$, $f(\sigma)\in L$. $K$ is simplicial isomorphic to $L$ if there are simplicial maps $f\colon K\to L$ and $g\colon L\to K$ such that $fg=1_L$ and $gf=1_K$ ($1_K$ and $1_L$ denote the identity map).
\vspace{3mm}

In the following definition, we do not use the definition of dual complex of a chain complex and regard the two complexes as the same chain group with different differentials $d$ and $\delta$. This technique is essential for the definition of $L$-spectral sequence and $L$-homology.\vspace{1.5mm}

{\bf Definition 2} Let $K$ be a simplicial complex with vertex set $S$.  $C_k(K)$ is the Abellian group generated by all ordered sequence $[v_0,v_1,\cdots,v_k]$ of elements of $S$ modular the following zero relations. $[v_0,\cdots,v_k]=0$ if $v_i=v_j$ for some $i\neq j$. If there is no repetition in $[v_0,\cdots,v_k]$, then $[v_0,\cdots,v_k]=0$ if $\{v_0,{\cdots},v_k\}\not\in K$ and $[\cdots,v_{i},\cdots,v_{j},\cdots]=-[\cdots,v_{j},\cdots,v_{i},\cdots]$ (the omitted part unchanged) for all $i<j$. Denote the unique generator of $C_{-1}(K)$ by $[\,]$. $[v_0,{\cdots},v_k]$ (including $[\,]$) is called a chain simplex of $K$ of dimension $k$. $C_*(K)=\oplus_{k\geqslant 0}C_k(K)$ and $\w C_*(K)=\oplus_{k\geqslant -1}C_k(K)$ are respectively the chain group and augmented chain group of $K$.

Two differentials $d\colon C_k(K)\to C_{k-1}(K)$ and $\delta\colon C_k(K)\to C_{k+1}(K)$ are defined as follows.
For any chain simplex $[v_0,\cdots,v_k]\in C_k(K)$ ($\hat v_i$ means canceling the symbol from the term),
\[d[v_0,\cdots,v_k]={\sum}_{i=0}^k(-1)^i[v_0,\cdots,\hat v_{i},\cdots,v_k]\,\,(d[v_0]=[\,]),\]
\[\delta[v_0,\cdots,v_k]={\sum}_{v\in S}[v,v_0,\cdots,v_n]\,\,(\delta([\,])={\sum}_{v\in S}[\,v\,]).\]

$H_*(K)=H_*(C_*(K),d)$ and $H^*(K)=H^*(C_*(K),\delta)$ are respectively the homology and cohomology of $K$; $\w H_*(K)=H_*(\w C_*(K),d)$ and $\w H^*(K)=H^*(\w C_*(K),\delta)$ are respectively the reduced homology and cohomology of $K$. For an Abellian group $G$, $H_*(K;G)=H_*(C_*(K){\otimes}G,d)$ and $H^*(K;G)=H^*(C_*(K){\otimes}G,\delta)$ are respectively the homology and cohomology of $K$ over $G$; $\w H_*(K;G)=H_*(\w C_*(K){\otimes}G,d)$ and $\w H^*(K;G)=H^*(\w C_*(K){\otimes}G,\delta)$ are respectively the reduced homology and cohomology of $K$ over $G$.
\vspace{3mm}

{\bf Convention} Notice the difference between simplex and chain simplex. For two simplexes $\sigma$ and $\tau$, $\sigma{\cap}\tau$ and $\sigma{\cup}\tau$ are naturally defined. But for two chain simplexes $\sigma$ and $\tau$, $\sigma{\cap}\tau$ and $\sigma{\cup}\tau$ are not naturally defined. We need only the definition of union of chain simplexes in later proofs. For two chain simplexes $\sigma=[v_0,{\cdots},v_m]$ and $\tau=[w_0,{\cdots},w_n]$, define $\sigma{\cup}\tau=[v_0,{\cdots},v_m,w_0,{\cdots},w_n]$ and $\tau{\cup}\sigma=[w_0,{\cdots},w_n,v_0,{\cdots},v_m]$. A chain simplex corresponds to a unique simplex consisting of all the vertices of the chain simplex. If there is no confusion, we use the same symbol to denote both the chain simplex and the simplex it corresponds to. So a chain simplex $\sigma\in K$ implies the corresponding simplex in $K$ and for two chain simplexes $\sigma$ and $\tau$, $\sigma\subset\tau$ and $\sigma{\cap}\tau=\phi$ implies the corresponding simplexes satisfy the condition. Since we only use the union $\sigma{\cup}\tau$ of chain simplexes when $\sigma{\cap}\tau=\phi$, $\sigma{\cup}\tau$ and $\tau{\cup}\sigma$ are chain simplexes that differ only up to sign.
\vspace{3mm}

{\bf Definition 3} For a simplicial complex $K$, the double complex $(T_{*,*}(K),\partial)$
\[\begin{array}{llllllll}&&&&&&T_{3,3}&\cdots\\&&&&&&_d\downarrow&\cdots\\
&&&&T_{2,2}&\stackrel{\delta}{\to}&T_{2,3}&\cdots\\&&&&_d\downarrow&&_d\downarrow&\cdots\\
&&T_{1,1}&\stackrel{\delta}{\to}&T_{1,2}&\stackrel{\delta}{\to}&T_{1,3}&\cdots\\
&&_d\downarrow&&_d\downarrow&&_d\downarrow&\cdots\\
T_{0,0}&\stackrel{\delta}{\to}&T_{0,1}&\stackrel{\delta}{\to}&T_{0,2}&\stackrel{\delta}{\to}&T_{0,3}&\cdots
\end{array}\]
is defined as follows. $ T_{*,*}(K)$ is the subgroup of $ C_*(K){\otimes} C_*(K)$ generated by all tensor product of chain simplexes $\sigma{\otimes}\tau$ such that $\sigma\!\subset\!\tau$ and $\partial\colon T_{s,t}\to T_{s-1,t}{\oplus}T_{s,t+1}$ is defined by $\partial(\sigma{\otimes}\tau)=(d\sigma){\otimes}\tau+(-1)^{|\sigma|}\sigma{\otimes}(\delta\tau)$. In another word, $(T_{*,*}(K),\partial)$ is the double subcomplex of $(C_*(K),d){\otimes}(C_*(K),\delta)$ generated by all tensor product of chain simplexes $\sigma{\otimes}\tau$ such that $\sigma{\subset}\tau$. The  $L$-spectral sequence $\{ L^r_{s,t}(K),\partial_r\}$ of $K$ is induced by the horizontal filtration $ L_{n}=\oplus_{s\leqslant n} T_{s,*}$ with $\partial_r\colon L^r_{s,t}\to L^r_{s-r,t-r+1}$. $ L^1_{*,*}(K)$ is simply denoted by $ L_{*,*}(K)$. $ L^2_{s,t}(K)$ denoted by $L H_{s,t}(K)$ is called the $L$-homology of $K$. For an Abellian group $G$, define $ T_{*,*}(K;G)= T_{*,*}(K){\otimes}G$ and we have similar definition of  $L$-spectral sequence $ L^r_{s,t}(K;G)$, $ L_{*,*}(K;G)= L^1_{*,*}(K;G)$ and $L$-homology $L H_{*,*}(K;G)=L^2_{*,*}(K;G)$ of $K$ over $G$.

Similarly, the augmented double complex $(\w T_{*,*}(K),\partial)$
\[\begin{array}{llllllll}&&&&&&\w T_{2,2}&\cdots\\&&&&&&_d\downarrow&\cdots\\
&&&&\w T_{1,1}&\stackrel{\delta}{\to}&\w T_{1,2}&\cdots\\&&&&_d\downarrow&&_d\downarrow&\cdots\\
&&\w T_{0,0}&\stackrel{\delta}{\to}&\w T_{0,1}&\stackrel{\delta}{\to}&\w T_{0,2}&\cdots\\
&&_d\downarrow&&_d\downarrow&&_d\downarrow&\cdots\\
\w T_{-1,-1}&\stackrel{\delta}{\to}&\w T_{-1,0}&\stackrel{\delta}{\to}&\w T_{-1,1}&\stackrel{\delta}{\to}&\w T_{-1,2}&\cdots
\end{array}\]
is the double subcomplex of $(\w C_*(K),d){\otimes}(\w C_*(K),\delta)$ generated by all product of chain simplexes $\sigma{\otimes}\tau$ such that $\sigma\!\subset\!\tau$. The reduced $L$-spectral sequence $\{\w L^r_{s,t}(K),\partial_r\}$ of $K$ is induced by the horizontal filtration $\w L_{n}=\oplus_{s\leqslant n}\w T_{s,*}$ with $\partial_r\colon\w L^r_{s,t}\to\w L^r_{s-r,t-r+1}$. $\w L^1_{*,*}(K)$ is simply denoted by $\w L_{*,*}(K)$. $\w L^2_{s,t}(K)$ denoted by $L\w H_{s,t}(K)$ is called the reduced $L$-homology of $K$. For an Abellian group $G$, define $\w T_{*,*}(K;G)=\w T_{*,*}(K){\otimes}G$ and we have similar definition of reduced $L$-spectral sequence $\w L^r_{s,t}(K;G)$, $\w L_{*,*}(K;G)=\w L^1_{*,*}(K;G)$ and reduced $L$-homology $L\w H_{*,*}(K;G)=\w L^2_{*,*}(K;G)$ of $K$ over $G$.
\vspace{3mm}

{\bf Theorem 4}\, Let $K$ be a simplicial complex and $G$ be an Abellian group.  The reduced $L$-spectral sequence of $K$ over $G$ satisfies that\\
\hspace*{15mm}$\w L_{s,s+t+1}(K;G)=\oplus_{|\sigma|=s}\, \w H^{t}({\rm link}_K\sigma;G)$,\,\,
$(\w L_{*,*}(K;G),\partial_1)=(\oplus_{\sigma\in K}\,\w H^{*}({\rm link}_K\sigma;G),\partial_1)$,\\
with $\partial_1$ on the right side defined as follows. Denote $x{\in}\w H^*({\rm link}_K\sigma;G)$ by $[x]_{\sigma}$ in $\oplus_{\sigma\in K}\,\w H^{*}({\rm link}_K\sigma;G)$, then $\partial_1([x]_{\sigma})=\sum_{i=0}^n[\psi_{\sigma}^{\sigma_i}(x)]_{\sigma_i}$, where $\sigma=\{v_0,{\cdots},v_n\}$, $\sigma_i=\sigma{-}\{v_i\}$ and $\psi_{\sigma}^{\sigma_i}$ is defined as follows. The group homomorphism $\w\psi^{\sigma_i}_{\sigma}\colon\w C_*({\rm link}_K\sigma;G)\to\w C_{*+1}({\rm link}_K\sigma_i;G)$ is defined by  $\w\psi_{\sigma}^{\sigma_i}(\tau)=\tau{\cup}[v_i]$ for all chain simplexes $\tau\in{\rm link}_K\sigma$. It is obvious that $\w\psi_{\sigma}^{\sigma_i}$ is a cochain complex homomorphism and $\psi_{\sigma}^{\sigma_i}=(\w\psi_{\sigma}^{\sigma_i})^*$ is the induced cohomology homomorphism.

The $L$-spectral sequence of $K$ over $G$ satisfies that \\
\hspace*{15mm}$L_{s,s+t+1}(K;G)=\oplus_{|\sigma|=s}\, \w H^{t}({\rm link}_K\sigma;G)$,\,\,
$(L_{*,*}(K;G),\partial_1)=(\oplus_{\sigma\in K,\sigma\neq\phi}\,\w H^{*}({\rm link}_K\sigma;G),\partial_1)$,\\
with $\partial_1$ on the right side defined as follows. For $x\in\w H^*({\rm link}_{K}\sigma;G)$, $\partial_1([x]_{\sigma})$ is just the same as the reduced case if $|\sigma|\!>\!0$ and $\partial_1([x]_{\sigma})=0$ if $|\sigma|\!=\!0$.

For any $n\!\geqslant\!0$, there is an exact sequence
\[0\to L\w H_{0,n}(K;G)\to LH_{0,n}(K;G)\to\w H^n(K;G)\to L\w H_{-1,n}(K;G)\to 0,\]
and $LH_{s,n}(K;G)=L\w H_{s,n}(K;G)$ if $s\!>\!0$.\vspace{3mm}

Proof\, $\w T_{u,v}(K)$ is generated by product of chain simplexes $\sigma\!\otimes\!\sigma'$ such that $\sigma\!\subset\!\sigma'$, $|\sigma|\!=\!u$, $|\sigma'|\!=\!v$. Equivalently, $\w T_{s,s+t+1}(K)$ is generated by all product of chain simplexes $\sigma\!\otimes\!(\tau{\cup}\sigma)$ such that $\sigma{\cap}\tau=\phi$, $|\sigma|=s$, $|\tau|=t$. So the correspondence $\sigma{\otimes}(\tau{\cup}\sigma)\to[\tau]_{\sigma}$ induces a graded group isomorphism $\varrho\colon\w T_{*,*}(K)\to \oplus_{\sigma\in K}\w C_*({\rm link}_K\sigma)$. Define $\partial$ on the latter group by  $\partial([\tau]_{\sigma})=
\sum_{i=1}^n[\w\psi^{\sigma_i}_{\sigma}(\tau)]_{\sigma_i}+(-1)^{|\sigma|}[\delta(\tau)]_{\sigma}$. Suppose $\sigma=[v_0,{\cdots},v_n]$ and $\tau=[w_0,{\cdots},w_m]$. Then
\begin{eqnarray*}&&\varrho\partial([v_0,{\cdots},v_n]{\otimes}[w_0,{\cdots},w_m,v_0,{\cdots},v_n])\\
&=&\varrho(\Sigma_{i=0}^n(-1)^{i}[v_0,\cdots,\hat v_i,\cdots,v_n]{\otimes}[w_0,{\cdots},w_m,v_0,{\cdots},v_n])\\
&&+\varrho(\Sigma_{\{w,w_0,\cdots,w_m\}\in{\rm link}_K\sigma}(-1)^{n}[v_0,\cdots,v_n]{\otimes}[w,w_0,\cdots,w_m,v_0,\cdots,v_n])\\
&=&\Sigma_{i=0}^n[w_0,\cdots,w_n,v_i]_{\{v_0,\cdots,\hat v_i,\cdots,v_n\}}+(-1)^n[\delta([w_0,\cdots,w_m])]_{\{v_0,\cdots,v_n\}}\\
&=&\partial([w_0,\cdots,w_m]_{\{v_0,\cdots,v_n\}})\\
&=&\partial\varrho([v_0,{\cdots},v_n]{\otimes}[w_0,{\cdots},w_m,v_0,{\cdots},v_n]).
\end{eqnarray*}
So $\varrho$ is a well-defined double complex isomorphism that induces a spectral sequence isomorphism with $E_2$-term isomorphism stated in the theorem.

Define trivial double complex $U_{*,*}(K;G)$ by $U_{-1,t}(K;G)=\w L_{-1,t}(K;G)=\w H^t(K;G)$ and $U_{s,t}(K;G)=0$ if $s\!\geqslant\!0$. Then $U_{*,*}$ is a subcomplex of $(\w L_{*,*},\partial_1)$ and the quotient complex is just $(L_{*,*},\partial_1)$. So from the long exact sequence induced by the short exact sequence $0\to U_{*,*}\to \w L_{*,*}\to L_{*,*}\to 0$ we get the relation between the reduced $L$-homology and $L$-homology. \vspace{3mm}

From the proof of the above theorem we know that there are two ways to describe the structure of $(T_{*,*}(K;G),\partial)$. One way is to denote the generators by $\sigma{\otimes}\sigma'$ with $\partial(\sigma{\otimes}\sigma')=(d\sigma){\otimes}\sigma'+(-1)^{|\sigma|}\sigma{\otimes}(\delta\sigma')$. The other is to denote the generators by $[\tau]_{\sigma}$ ($\tau\in{\rm link}_K\sigma$) with $\partial([\tau]_{\sigma})=\sum_{i=1}^n[\w\psi^{\sigma_i}_{\sigma}(\tau)]_{\sigma_i}+
(-1)^{|\sigma|}[\delta(\tau)]_{\sigma}$. We have to use both in later proofs. But we only denote the generators of $(L_{*,*}(K;G),\partial_1)$ in one way, i.e., $[x]_{\sigma}$ with $x\in\w H^*({\rm link}_K\sigma;G)$ and $\partial_1([x]_{\sigma})=\sum_{i=0}^n[\psi^{\sigma_i}_{\sigma}(x)]_{\sigma_i}
=\sum_{i=0}^n[x{\cup}[v_i]]_{\sigma_i}$. The reduced case is the same.\vspace{3mm}

{\bf Remark} We can similarly define $U_{s,t}(K)=C_s(K){\otimes}C_t(K^c)$ generated by all $\sigma{\otimes}\tau$ with $\sigma\in K$ and $\tau\not\in K$ such that $\sigma\subset\tau$ and $\partial(\sigma{\otimes}\tau)=(d\sigma){\otimes}\tau+(-1)^{|\sigma|}\sigma{\otimes}(\delta\tau)$. Then the $R$-spectral sequence of $K$ $\{R_r^{s,t}(K),\partial^r\}$ induced by the vertical filtration $R^{n}=\oplus_{s\geqslant n}U_{*,s}$ satisfies $R_1^{*,*}(K)=\oplus_{\sigma\not\in K}\w H_*(K|_{\sigma})$, where $K|_{\sigma}=\{\tau{\in}K\,|\,\tau{\subset}\sigma\}$. But $R$-homology $RH^{s,t}(K;G)=R_2^{s,t}(K;G)$ is neither a homeomorphism or homotopy invariant of $K$  nor a homeomorphism or homotopy invariant of the Alexander dual $K^*$ of $K$.\vspace{3mm}

{\bf Definition 5} For an Abellian group $G$, a simplicial complex $K$ is of essential cohomology dimension $n$ ($n\!>\!0$) over $G$ if for all non-empty simplex $\sigma$ of $K$, $\w H^t({\rm link}_K\sigma;G)=0$ for all $t\neq n{-}|\sigma|{-}1$.\vspace{3mm}

{\bf Theorem 6} Let $K$ be a simplicial complex of essential cohomology dimension $n$ over $G$. Then for all $0\!\leqslant\!s\!<\!n$, $L\w H_{s,n}(K;G)=L\w H_{-1,n-s}(K;G)=\w H^{n-s}(K;G)$, $L\w H_{n,n}(K;G)=G$ and $L\w H_{s,t}(K;G)=0$ otherwise.\vspace{5mm}

Proof\, By definition, $\w L_{-1,s}(K;G)=\w H^s(K;G)$ for all $0\!\leqslant\!s\!\leqslant\!n$ and $\w L_{s,t}(K;G)=0$ if $s\!\neq\!-1$ or $t\!\neq\!n$. So the reduced $L$-spectral sequence satisfies $\w L^{s+1}_{s,n}(K;G)=\w L_{s,n}(K;G)$ and $\w L^{s+1}_{-1,n-s}(K;G)=\w L^{s+1}_{-1,n-s}(K;G)$ for all $0\!<\!s\!\leqslant\!n$. From the vertical filtration $E_k=\oplus_{t\geqslant k}T_{*,t}(K;G)$ we get a spectral sequence $E^r_{s,t,k}$ with $E^2_{*,*,*}=\oplus_{\sigma\in K}\w H^*(2^{\sigma};G)=\w H^*(2^{\phi};G)=G$ ($2^{\sigma}$ as defined in the following Example 7). So $H^*(\w T_{*,*}(K;G),\partial)=G$ with the generator class represented by $\sum_{\sigma\in K}(-1)^{[\frac{|\sigma|}{2}]}\sigma{\otimes}\sigma$.  So the reduced $L$-spectral sequence converges to $G$ and $\partial_{s+1}\colon\w L_{s,n}(K;G)\to\w L_{-1,n-s}(K;G)$ is an isomorphism for all $0\!\leqslant\!s\!<\!n$ and $\partial_{n+1}\colon\w L_{n,n}(K;G)\to\w L_{-1,0}(K;G)$ is an epimorphism with kernel $G$.\vspace{3mm}

{\bf Example 7} For a set $\sigma=\{v_0,\cdots,v_n\}$, $2^{\sigma}$ is the simplicial complex consisting of all subsets of $\sigma$ and $\partial\sigma$ is the simplicial complex consisting of all proper subsets of $\sigma$. By Theorem 4, $L\w H_{n,n}(2^{\sigma};G)=\w H^{-1}({\rm link}_{2^{\sigma}}\sigma;G)=\w H^{-1}(\{\phi\};G)=G$ with generator class represented by $[\phi]_{\sigma}$ and $L\w H_{s,t}(2^{\sigma};G)=0$ otherwise. For all $\tau\in\partial\sigma$, ${\rm link}_{\partial\sigma}\tau=\partial(\sigma{-}\tau)$. So for $n\!\geqslant\!2$, $\partial\sigma$ is of essential cohomology dimension $n{-}1$. By Theorem 6, $L\w H_{n-1,n-1}(\partial\sigma;G)=G$ and $L\w H_{s,t}(\partial\sigma;G)=0$ otherwise. A direct checking shows that the reduced $L$-homology formula holds for all $n\!\geqslant\!0$ and the generator class of $L\w H_{n-1,n-1}(\partial\sigma;G)$ is represented in $\w L_{n-1,n-1}(\partial\sigma;G)$ by $\sum_{i=0}^n[\phi]_{\sigma_i}$ with $\sigma_i=\sigma{-}\{v_i\}$.\vspace{3mm}

{\bf Definition 8} The geometrical realization of a simplicial complex $K$ is the topological space $|K|$ defined as follows. Suppose the vertex set $S=\{v_1,\cdots,v_m\}$. Let $e_i\in\RR^m$ be such that the $i$-th coordinate of $e_i$ is $1$ and all other coordinates of $e_i$ are $0$. Then $|K|$ is the union of all the convex hull of $e_{i_0},\cdots,e_{i_s}$ such that $\{v_{i_0},\cdots,v_{i_s}\}$ is a simplex of $K$. The convex hull of $e_{i_0},\cdots,e_{i_s}$ is a geometrical $s$-simplex of $|K|$. A polyhedron $X$ is a topological space that is homeomorphic to the geometrical realization $|K|$ of a simplicial complex $K$ and $K$ is a triangulation of $X$.\vspace{3mm}

From Example 7 we know that $LH_{*,*}(K;G)$ and $L\w H_{*,*}(K;G)$ are not homotopy invariant of the polyhedron $|K|$ since $|2^{\sigma}|$ is contractible but $L\w H_{*,*}(2^{\sigma};G)\neq L\w H_{*,*}(2^{\tau};G)$ when $|\sigma|\!\neq\!|\tau|$. In fact, simplicial maps do not induce (reduced) $L$-homology homomorphisms. However, $LH_{*,*}$ and $L\w H_{*,*}$ are homeomorphism invariant of polyhedra.\vspace{3mm}

{\bf Definition 9} Let $K$ be a simplicial complex. A stellar subdivision of $K$ on the simplex $\sigma$ is the simplicial complex $S_{\sigma}(K)$ defined as follows. If $\sigma=\phi$, then $S_{\phi}(K)=K$. If $|\sigma|=0$, then $S_{\sigma}(K)$ is simplicial isomorphic to $K$ by replacing the unique vertex $w$ of $\sigma$ by a new vertex $v$. If $|\sigma|>0$, then the vertex set of $S_{\sigma}(K)$ is the vertex set of $K$ added with a new vertex $v$ and ($\partial\sigma$ is as defined in Example 7)\\
\hspace*{12.5mm}$S_{\sigma}(K)=(K{-}\{\sigma\}{*}{\rm link}_K\sigma)\cup C(\partial\sigma{*}{\rm link}_K\sigma)=\{\tau\!\in\!K\,|\,\sigma\!\not\subset\!\tau\}\cup
\{\{v\}{\cup}\sigma'{\cup}\tau\,|\,\sigma'\!\in\!\partial\sigma,\,\tau\!\in\!{\rm link}_K\sigma
\}$.\vspace{3mm}

Suppose the vertex set of $K$ is $\{v_1,\cdots,v_m\}$ and  $\sigma=\{v_{i_0},\cdots,v_{i_s}\}\in K$. Then the map $f(e_i)=e_i$ for $i=1,\cdots,m$, $f(e_{m+1})=\frac {1}{s+1}(e_{i_0}{+}\cdots{+}e_{i_s})$ extends linearly to a homeomorphism from $|S_{\sigma}(K)|$ to $|K|$. Conversely, by \cite{a} and \cite{b}, two polyhedra $|K|$ and $|L|$ are homeomorphic if and only if there are simplicial complexes $K_0=K,K_1,\cdots,K_n,K_{n+1}=L$ such that for every $i=0,\cdots,n$, either $K_{i+1}$ is a stellar subdivision of $K_i$, or $K_{i}$ is a stellar subdivision of $K_{i+1}$. \vspace{3mm}

{\bf Theorem 10} Let $K$ be a simplicial complex and $G$ be an Abellian group. Then for any simplex $\sigma$ of $K$, $LH_{*,*}(K;G)=LH_{*,*}(S_{\sigma}(K);G)$ and $L\w H_{*,*}(K;G)=L\w H_{*,*}(S_{\sigma}(K);G)$. Thus for a polyhedron $X$, define $LH_{*,*}(X;G)=LH_{*,*}(K;G)$ and $L\w H_{*,*}(X;G)=L\w H_{*,*}(K;G)$ for any triangulation $K$ of $X$. Then $LH_{*,*}(X;G)$ and $L\w H_{*,*}(X;G)$ are independent of the triangulation and are homeomorphism invariants of $X$.\vspace{1.5mm}

Proof\, We first prove the reduced case. If $|\sigma|\leqslant 0$, the conclusion holds by definition. Suppose $\sigma=\{v_0,\cdots,v_n\}$ ($=[v_0,\cdots,v_n]$ when necessary) with $n\!>\!0$. Denote $L=S_{\sigma}(K)$.

Let $M=K{\cap}L$ be the simplicial subcomplex of $K$ and $L$ consisting of all simplexes $\tau$ such that $\sigma{\cap}\tau$ is a proper face of $\sigma$. Let $\w T_{*,*}(K|_M)$ and $\w T_{*,*}(L|_M)$ be respectively the subgroup of $\w T_{*,*}(K)$ and $\w T_{*,*}(L)$ generated by all tensor product of chain simplexes $\sigma{\otimes}\tau$ such that $\sigma\!\in\!M$ and $\sigma\!\subset\!\tau$. Then $(\w T_{*,*}(K|_M),\partial)$ and $(\w T_{*,*}(L|_M),\partial)$ are respectively double subcomplexes of $(\w T_{*,*}(K),\partial)$ and $(\w T_{*,*}(L),\partial)$ generated by all tensor product of chain simplexes $\sigma{\otimes}\tau$ such that $
\sigma\in M$. From the horizontal filtration we get reduced $L$-spectral sequences $\w L^r_{s,t}(K|_M;G)$ and $\w L^r_{s,t}(L|_M;G)$ and reduced $L$-homology $L\w H_{s,t}(K|_M;G)$ and $L\w H_{s,t}(L|_M;G)$. Define $(\w T_{*,*}(K/M),\partial)=(\w T_{*,*}(K),\partial)/(\w T_{*,*}(K|_M),\partial)$ and $(\w T_{*,*}(L/M),\partial)=(\w T_{*,*}(L),\partial)/(\w T_{*,*}(L|_M),\partial)$. From the horizontal filtration we also get reduced $L$-spectral sequences $\w L^r_{s,t}(K/M;G)$ and $\w L^r_{s,t}(L/M;G)$ and reduced $L$-homology $L\w H_{s,t}(K/M;G)$ and $L\w H_{s,t}(L/M;G)$. We have the following two short exact sequences of double complexes\\
\hspace*{20mm}$\begin{array}{ccccccccc}
    0&\to & \w L_{*,*}(K|_M;G) & \to &\w L_{*,*}(K;G)& \to & \w L_{*,*}(K/M;G) & \to &0 \\
    0&\to & \w L_{*,*}(L|_M;G) & \to &\w L_{*,*}(L;G)& \to & \w L_{*,*}(L/M;G) & \to &0
  \end{array}$\\

For all simplex $\tau\!\in\!M$, ${\rm link}_L\tau=S_{\sigma-\tau}({\rm link}_K\tau)$ and let $f_{\tau}\colon {\rm link}_L\tau\to{\rm link}_K\tau$ be the simplicial map defined by $f_{\tau}(v)=v_s$ for a fixed $v_s\in\sigma{-}\tau$ and $f_{\tau}(w)=w$ for all other $w$ ($f_{\tau}$ is the identity map if $\sigma{-}\tau=\phi$). Then $f_{\tau}$ induces a homomorphism $\w f_{\tau}\colon(\w C_*({\rm link}_K\tau;G),\delta)\to(\w C_*({\rm link}_L\tau;G),\delta)$ and an isomorphism $f^*_{\tau}\colon H^*({\rm link}_K\tau)\to H^*({\rm link}_L\tau)$.
Theorem 4 also holds for $K|_M$ and $L|_M$. Precisely, we have two isomorphisms\\
\hspace*{10mm}$(\w T_{*,*}(K|_M;G),\partial)=(\oplus_{\tau\in M}\w C_*({\rm link}_K\tau;G),\partial)\quad\quad
(\w T_{*,*}(L|_M;G),\partial)=(\oplus_{\tau\in M}\w C_*({\rm link}_L\tau;G),\partial)$\\
with $\partial([\eta]_{\tau})=(-1)^{\tau}[\delta(\eta)]_{\tau}+\sum_i[\eta{\cup}[w_i]]_{\tau_i}$ for all chain simplex $\eta\in {\rm link}_K\tau$ or ${\rm link}_L\tau$ and $\tau_i=\tau{-}\{w_i\}$. $\varpi'=\oplus_{\tau\in M}\w f_{\tau}$ is a graded group isomorphism between the two double complexes. To prove $\varpi'\partial=\partial\varpi'$, we may suppose $\tau=\{v_s,{\cdots},v_n,v_{n+1},{\cdots},v_{n+t}\}$ and $f_{\tau}(v)=v_1$. Denote $\tau_i=\tau{-}\{v_i\}$. Then for a chain simplex $\eta\in{\rm link}_K\tau$ such that $\sigma{-}\tau$ is not a proper face of $\eta$, $\w f_{\tau}(\eta)=\eta$ and $\w\psi_{\tau}^{\tau_i}\w f_{\tau}(\eta)=\w f_{\tau}\w\psi_{\tau}^{\tau_i}(\eta)=\w\psi_{\tau}^{\tau_i}(\eta)$. For $\eta\in{\rm link}_K\tau$ such that $\eta=[v_1,{\cdots},v_{s-1},v_{l_1},{\cdots},v_{l_u}]$ ($l_j\!>\!n{+}t$), $\w\psi_{\tau}^{\tau_i}\w f_{\tau}(\eta)=\w\psi_{\tau}^{\tau_i}([v,v_2,{\cdots},v_{s-1},v_{l_1},{\cdots},v_{l_u}])
=[v,v_2,{\cdots},v_{s-1},v_{l_1},{\cdots},v_{l_u},v_i]=\w f_{\tau}([v_1,{\cdots},v_{s-1},v_{l_1},{\cdots},v_{l_u},v_i])=\w f_{\tau}\w\psi_{\tau}^{\tau_i}(\eta)$.
Thus $\w\psi_{\tau}^{\tau_i}\w f_{\tau}=\w f_{\tau}\w\psi_{\tau}^{\tau_i}$. Since $\w f_{\tau}$ is induced by simplicial map, $\delta\w f_{\tau}=\w f_{\tau}\delta$. So $\varpi'\partial=\partial\varpi'$ and $\varpi'$ is a double complex homomorphism from  $(\w T_{*,*}(K|_M;G),\partial)$ to $(\w T_{*,*}(L|_M;G),\partial)$. Since $f^*_{\tau}$ is an isomorphism for all $\tau\in M$, $\varpi_1=\oplus_{\tau\!\in\!M}f^*_{\tau}$ is a double complex isomorphism from $(\w L_{*,*}(K|_M;G),\partial_1)$ to $(\w L_{*,*}(L|_M;G),\partial_1)$.

Denote $B={\rm link}_K\sigma$. Then for $\tau\in B$ and $\sigma'\in\partial\sigma$ ($\partial\sigma$ as defined in Example 7),\\
\hspace*{30mm}${\rm link}_K\sigma{\cup}\tau={\rm link}_B\tau,\quad\quad{\rm link}_L\{v\}{\cup}\sigma'{\cup}\tau=({\rm link}_{\partial\sigma}\sigma')*({\rm link}_B\tau)$.\\
The two simplicial isomorphisms implies two isomorphisms $\varrho_1\colon (\w T_{n+*,n+*}(K/M),\partial)\to (\w T_{*,*}(B),\partial)$ and $\varrho_2\colon(\w T_{*+1,*+1}(L/M),\partial)\to (\w T_{*,*}(\partial\sigma),\partial)\otimes(\w T_{*,*}(B),\partial)$ defined as follows.
$\w T_{*,*}(K/M)$ is generated by all tensor product of chain simplexes $\tau{\cup}\sigma\otimes\eta{\cup}\tau{\cup}\sigma$ such that $\tau\in B$ and $\eta\in{\rm link}_B\tau$. So the correspondence $\tau{\cup}\sigma\otimes\eta{\cup}\tau{\cup}\sigma\to (-1)^{|\sigma|+1}\tau\otimes\eta{\cup}\tau$ induces $\varrho_1$. Similarly, $\w T_{*,*}(L/M)$ is generated by all tensor product of chain simplexes $\sigma'{\cup}\tau{\cup}[v]\otimes\xi{\cup}\eta{\cup}\sigma'{\cup}\tau{\cup}[v]$ such that $\xi\in{\rm link}_{\partial\sigma}\sigma'$ and $\eta\in {\rm link}_B\tau$. So the correspondence $\sigma'{\cup}\tau{\cup}[v]\otimes\xi{\cup}\eta{\cup}\sigma'{\cup}\tau{\cup}[v]\,\to\, (-1)^{(|\xi|+1)|\tau|+(|\eta|+1)|\sigma'|}(\sigma'{\otimes}\xi{\cup}\sigma')\otimes(\tau{\otimes}\eta{\cup}\tau)$ induces $\varrho_2$. In another way, $\varrho_1([\eta]_{\sigma{\cup}\tau})=(-1)^{|\sigma|+1}[\eta]_{\tau}$ for all chain simplex $\eta\in{\rm link}_K\sigma{\cup}\tau\!=\!{\rm link}_B\tau$ and $\varrho_2([\xi{\cup}\eta]_{\{v\}\cup\sigma'\cup\tau})=(-1)^{(|\xi|+1)|\tau|+(|\eta|+1)|\sigma'|}
[\xi]_{\sigma'}{\otimes}[\eta]_{\tau}$ for all $\xi{\cup}\eta\in{\rm link}_L\{v\}{\cup}\sigma'{\cup}\tau\!=\!{\rm link}_{\partial\sigma}\sigma'*{\rm link}_B\tau$. Since $L\w H_{*,*}(\partial\sigma)=\zz$ represented by $\sum_{i=0}^n[\phi]_{\sigma_i}$ (see Example 7), the correspondence\\
\hspace*{60mm}$[x]_{\sigma{\cup}\tau}\,\to\,
\sum_{i=0}^n[x]_{\{v\}{\cup}\sigma_i{\cup}\tau}$\\
for all $x\in\w H^*({\rm link}_B\tau;G)$ and $\tau\in B$ induces a double complex homomorphism $\varpi_2\colon(\w L_{*,*}(K/M;G),\partial_1)\to(\w L_{*,*}(L/M;G),\partial_1)$ that induces a reduced  $L$-homology isomorphism.

Since there are graded group isomorphism $\w L_{*,*}(K;G)=\w L_{*,*}(K|_M;G){\oplus}\w L_{*,*}(K/M;G)$ and $\w L_{*,*}(L;G)=\w L_{*,*}(L|_M;G){\oplus}\w L_{*,*}(L/M;G)$, we have a graded group homomorphism $\varpi=\varpi_1{\oplus}\varpi_2\colon\w L_{*,*}(K;G)\to\w L_{*,*}(L;G)$. By the previous conclusion, $\partial_1\varpi(z)=\varpi\partial_1(z)$ for all $z\!\in\!\w L_{*,*}(K|_M;G)$. Notice that $\w L_{*,*}(K/M;G)$ is generated by element of the form $[x]_{\sigma\cup\tau}$ with $x\in\w H^*({\rm link}_K\sigma{\cup}\tau;G)\!=\!\w H^*({\rm link}_B\tau;G)$. Suppose $\tau=\{w_1,{\cdots},w_m\}$. Then
\begin{eqnarray*}&&\partial_1\varpi([x]_{\sigma\cup\tau})\\
&=&\partial_1(\Sigma_{i}\,[x]_{[v]\cup\sigma_i\cup\tau})\\
&=&\Sigma_{i}\,[x{\cup}[v]]_{\sigma_i\cup\tau}
+\Sigma_{i<j}\,[x{\cup}([v_i]{+}[v_j])]_{\{v\}\cup\sigma_{i,j}\cup\tau}+
\Sigma_{i,k}\,[x{\cup}[w_k]]_{\{v\}{\cup}\sigma_i{\cup}(\tau-\{w_k\})}\\
&=&\Sigma_{i}\,[x{\cup}[v]]_{\sigma_i\cup\tau}
+\Sigma_{i,k}\,[x{\cup}[w_k]]_{\{v\}{\cup}\sigma_i{\cup}(\tau-\{w_k\})}\\
&=&\Sigma_{i}\,\varpi_1([x{\cup}[v_i]]_{\sigma_i\cup\tau})
+\Sigma_{k}\,\varpi_2([x{\cup}[w_k]]_{\sigma{\cup}(\tau-\{w_k\})})\\
&=&\varpi\partial_1([x]_{\sigma\cup\tau})
\end{eqnarray*}
where $\sigma_{i,j}=\{v_0,{\cdots},\hat v_i,{\cdots},\hat v_j,{\cdots},v_n\}$ and $\Sigma_{i<j}\,[x{\cup}(\{v_i\}{+}\{v_j\})]_{\{v\}\cup\sigma_{i,j}\cup\tau}
=\sum_i(-1)^{|x|+1}[\delta(x)]_{\{v\}{\cup}\sigma_{i,j}{\cup}\tau}=0$. Thus, $\partial_1\varpi=\varpi\partial_1$ and we have the following commutative diagram of short exact sequences of double complexes
\begin{eqnarray*}\begin{array}{ccccccccc}
    0&\to & \w L_{*,*}(K|_M;G) & \to &\w L_{*,*}(K;G)& \to & \w L_{*,*}(K/M;G) & \to &0 \\
    &&\varpi_1\downarrow\quad&&\varpi\downarrow\quad&&\varpi_2\downarrow\quad&&\\
    0&\to & \w L_{*,*}(L|_M;G) & \to &\w L_{*,*}(L;G)& \to & \w L_{*,*}(L/M;G) & \to &0
  \end{array}\end{eqnarray*}
that induces a commutative diagram of long exact sequences of reduced $L$-homology
\[\begin{array}{ccccccccl}\cdots\to& L\w H_{s,t}(K|_M;G)&{\to}&L\w H_{s,t}(K;G)&{\to}&L\w H_{s,t}(K{/}M;G)&
{\to}&L\w H_{s-1,t}(K|_M;G)&\to\cdots\\
&\varpi_1^*\downarrow\quad &&\varpi^*\downarrow\,\,&&\varpi_2^*\downarrow\quad &
&\,\,\varpi_1^*\downarrow\quad &\\
\cdots\to&L\w H_{s,t}(L|_M;G)&{\to}&L\w H_{s,t}(L;G)&{\to}&L\w H_{s,t}(L{/}M;G)&
{\to}&L\w H_{s-1,t}(L|_M;G)&\to\cdots.
\end{array}\]
Since $\varpi_1^*$ and $\varpi_2^*$ are isomorphisms, $\varpi^*\colon L\w H_{*,*}(K;G)\to L\w H_{*,*}(L;G)$ is an isomorphism.

Let $U_{*,*}$ be as in the proof of Theorem 4. Then  $U_{-1,*}(K;G)=\w H^*(K;G)=\w H^*(L;G)=U_{-1,*}(L;G)$, so the reduced $L$-homology isomorphism  from $K$ to $L$ induces $L$-homology isomorphism.
\vspace{3mm}

By Theorem 10, the disk $D^m$ and $D^n$ are not homeomorphic if $m\neq n$, since their (reduced) $L$-homology are different.\vspace{3mm}

{\bf Theorem 11} For two polyhedra $X$ and $Y$ and Abellian group $G$,
\[\begin{array}{ccc}
    LH_{*,*}(X\sqcup Y;G)&=&LH_{*,*}(X;G)\oplus LH_{*,*}(Y;G), \\
    L\w H_{*,*}(X\vee Y;G)&=&L\w H_{*,*}(X;G)\w\oplus L\w H_{*,*}(Y;G),
  \end{array}\]
where $\sqcup$  and $\vee$ are respectively the disjoint union and one point union of topological spaces and $\w\oplus$ is defined by $L\w H_{-1,0}(X\vee Y;G)=G{\oplus}L\w H_{-1,0}(X;G)\oplus L\w H_{-1,0}(Y;G)$ and $L\w H_{s,t}(X\vee Y;G)=L\w H_{s,t}(X;G)\oplus L\w H_{s,t}(Y;G)$ otherwise.
\vspace{1.5mm}

Proof\, Since for any triangulation $K$ of $X$ and $L$ of $Y$, we have by Theorem 4\\
\hspace*{36mm}$(T_{*,*}(K\sqcup L;G),\partial\,)=(T_{*,*}(K;G),\partial\,){\oplus}(T_{*,*}(L;G),\partial\,)$,\\
\hspace*{35mm}$(\w L_{*,*}(K\vee L;G),\partial_1)=(\w L_{*,*}(K;G),\partial_1){\w\oplus}(\w L_{*,*}(L;G),\partial_1)$,\\
the theorem holds.\vspace{3mm}

{\bf Theorem 12} For two polyhedra $X$ and $Y$ and a field $\ff$,
\[\begin{array}{ccc}
    L\w H_{*,*}(X{*}Y;\ff)&=&L\w H_{*,*}(X;\ff)\,{\w\otimes}\,L\w H_{*,*}(Y;\ff), \\
        LH_{*,*}(X{\times} Y;\ff)&=&LH_{*,*}(X;\ff)\,{\otimes}\, LH_{*,*}(Y;\ff),
  \end{array}\]
where $\w\otimes$ is defined by $L\w H_{s,t}(X{*}Y)=\oplus_{s'+s''=s-1,\,t'+t''=t-1}L\w H_{s',t'}(X){\otimes}L\w H_{s'',t''}(Y)$. The cone satisfies that for all Abellian group $G$, $L\w H_{s+1,t+1}(CX;G)=L\w H_{s,t}(X;G)$ for all $s,t\geqslant 0$ and $L\w H_{u,v}(CX;G)=0$ if $u<1$ or $v<1$. \vspace{1.5mm}

Proof\, For any two simplicial complexes $K$ and $L$, $(\w T_{*,*}(K*L),\partial)=(\w T_{*,*}(K),\partial){\w\otimes}(\w T_{*,*}(L),\partial)$. For any two simplexes $\sigma\in K$ and $\tau\in L$, ${\rm link}_{K*L}\sigma{\sqcup}\tau={\rm link}_K\sigma*{\rm link}_L\tau$. So the union case holds by Theorem 4 and K\"{u}nneth Theorem. If $K=\{\phi,\{v\}\}$, the cone case holds for all Abellian group $G$ by definition.

The Cartesian product case $K\bar\times L$ is complicated and we only give a sketch of the proof. Suppose the ordered vertex sets of $K$ and $L$ are respectively $S=\{v_1,\cdots,v_m\}$ and $T=\{w_1,\cdots,w_n\}$. Denote $\o{S{\times}T}=S{\times}T\cup\{\phi\}$. Define relation $\prec$ on $\o{S{\times}T}$ as follows. $(v_i,w_j)\prec(v_k,w_l)$ if $i\!\leqslant\!k$ and $j\!\leqslant\!l$ and $(i,j)\!\neq\!(k,l)$. $\phi\prec (v_i,w_j)$ and $(v_i,w_j)\prec\phi$ for all $(v_i,w_j)$. For $b_1,b_2\!\in\!\o{S{\times}T}$ with $b_1=(v_i,w_j)$ and $b_2=(v_k,w_l)$ ($i,j=0$ if $b_1=\phi$ and $k,l=m{+}n{+}1$ if $b_2=\phi$),
\begin{eqnarray*}B(b_1,b_2)\!\!\!&=&\!\!\!\{\xi\!\in\! K{\bar\times}L\,|\,\xi=\{\phi\},\,{\rm or}\,\,\xi\!=\!\{c_1,{\cdots},c_n\}\,{\rm such\,that}\,\,b_1\!\prec\!c_1\!\prec\!\cdots\!\prec\!c_n\!\prec\!b_2\},\\
B_K(b_1,b_2)\!\!\!&=&\!\!\!\{\,\,\,\sigma\in K\,\,\,|\,\sigma=\{\phi\},\,{\rm or}\,\,\sigma\!=\!\{v_{i_1},{\cdots},v_{i_s}\}\,{\rm such\,that}\,\,i\!<\!i_1\!<\!\cdots\!<\!i_s\!<\!k\},\\
B_L(b_1,b_2)\!\!\!&=&\!\!\!\{\,\,\,\tau\in\, L\,\,\,\,|\,\tau=\{\phi\},\,{\rm or}\,\,\tau\!=\!\{v_{j_1},{\cdots},v_{j_t}\}\,{\rm such\,that}\,\,j\!<\!j_1\!<\!\cdots\!<\!j_t\!<\!l\}.\end{eqnarray*}
Then their geometrical realizations satisfy
\begin{eqnarray*}|B(b_1,b_2)|\!\!\!&\cong&\!\!\!|B_K(b_1,b_2))|\barwedge|B_L(b_1,b_2)|\,\,{\rm if\,one\, of}\,\,b_1\,\,{\rm and}\,\,b_2\,\,{\rm is}\,\,\phi,\\
|B(b_1,b_2)|\!\!\!&\cong&\!\!\!|B_K(b_1,b_2))|\doublebarwedge|B_L(b_1,b_2)|\,\,{\rm if}\,\,b_1,b_2\neq\phi\,\,{\rm and}\,\,i<k,\,j<l,\\
|B(b_1,b_2)|\!\!\!&\cong&\!\!\!|B_K(b_1,b_2))|\,\,{\rm if}\,\,b_1,b_2\neq\phi\,\,{\rm and}\,\,i<k,j=l,\\
|B(b_1,b_2)|\!\!\!&\cong&\!\!\!|B_L\,(b_1,b_2))|\,\,{\rm if}\,\,b_1,b_2\neq\phi\,\,{\rm and}\,\,i=k,j<l,\end{eqnarray*}
where
\begin{eqnarray*}X\!\barwedge Y\!\!&=&\!\!(CX{\times}Y\cup X{\times}CY)/\!\!\sim{\rm with}\,\,((1,x),y)\!\sim\!(x,(1,y))\,{\rm for\, all}\,(x,y)\!\in\!X{\times}Y,\\
X\!\doublebarwedge Y\!\!&=&\!\!(CX{\barwedge}\,Y\cup X{\barwedge}\,CY)/\!\!\sim\,{\rm with}\,(1,x)\barwedge \,y\!\sim x\barwedge(1,y)\,{\rm for\, all}\,\,x\barwedge y\in\!X\barwedge Y,\end{eqnarray*}
and $CX=I{\times}X/\{0\}{\times}X$, $CY=I{\times}Y/\{0\}{\times}Y$. From the Mayer-Vietorus sequence we have that
\begin{eqnarray*}
\w H^*(B(b_1,b_2);\ff)\!\!\!&=&\!\!\!\Sigma^2\w H^*(B_K(b_1,b_2);\ff)\otimes\w H^*(B_L(b_1,b_2);\ff)\,\,{\rm if}\,\,b_1,b_2\neq\phi\,\,{\rm and}\,\,i<k,\,j<l,\\
\w H^*(B(b_1,b_2);\ff)\!\!\!&=&\!\!\!\Sigma\,\,\,\w H^*(B_K(b_1,b_2);\ff)\otimes \w H^*(B_L(b_1,b_2);\ff)\,\,{\rm otherwise},\end{eqnarray*}
where $\Sigma$ means uplifting the degree by 1.

For a non-empty simplex $\xi=\{b_0,{\cdots},b_n\}\in K{\bar\times}L$ such that $b_0\!\prec\!{\cdots}\!\prec\!b_n$,\\
\hspace*{31mm}${\rm link}_{K\bar\times L}\xi=B(\phi,b_0)*B(b_0,b_1)*\cdots*B(b_{n-1},b_n)*B(b_n,\phi).$\\
Define $P(\xi)=(\sigma,\tau)\in K{\times}L$ to be such that $\xi\subset\sigma{\times}\tau$ and $\xi\not\subset\sigma'{\times}\tau$, $\xi\not\subset\sigma{\times}\tau'$ for all proper faces $\sigma'$ of $\sigma$ and $\tau'$ of $\tau$. Then we have\\
\hspace*{30mm}$\w H^*({\rm link}_{K\bar\times L}\xi;\ff)=\Sigma^{|\sigma|+|\tau|-|\xi|+1}\w H^*({\rm link}_{K}\sigma;\ff)\otimes\w H^*({\rm link}_{L}\tau;\ff)$\\
and a direct checking shows that the homomorphism $\psi^{\xi_i}_{\xi}$ with $\xi_i=\xi{-}\{b_i\}$ is an isomorphism if $P(\xi_i)=P(\xi)$.

Thus, by Theorem 4, there is a double complex isomorphism\\
\hspace*{15mm}$(L_{*,*}(K{\bar\times}L;\ff),\partial_1)=(\oplus_{P(\xi)=
(\sigma,\tau),\xi\neq\phi}\Sigma^{|\sigma|+|\tau|-|\xi|+1}\w H^*({\rm link}_K\sigma;\ff){\otimes}\w H^*({\rm link}_L\tau;\ff),\w\partial_1)$.\\
Give the right side double complex a filtration $\{E_n\}$ as follows. For $[x]_{\xi}$ such that $P(\xi)=(\sigma,\tau)$, define $\|[x]_{\xi}\|=|\sigma|{+}|\tau|$ and $E_n$ to be the subgroup generated by all $[x]_{\xi}$ such that $\|[x]_{\xi}\|\leqslant n$. Then we get a spectral sequence $(E_{s,t,u}^r,d_r)$ from the filtration converging to $LH_{*,*}(K{\bar\times}L;\ff)$ such that\\
\hspace*{6mm}$(E^1_{*,*,*},d_1)=\oplus_{(\sigma,\tau)\in K{\times}L,\sigma,\tau\neq\phi}(C_*(\sigma{\bar\times}\tau)/C_*(\partial(\sigma{\bar\times}\tau)),d)
{\otimes}\Sigma^{|\sigma|+|\tau|-|\xi|+1}\w H^*({\rm link}_K\sigma;\ff){\otimes}\w H^*({\rm link}_L\tau;\ff)$.\\
Since $H_*(C_*(\sigma{\bar\times}\tau)/C_*(\partial(\sigma{\bar\times}\tau));\ff)=\ff$ with generator class represented by
$\Sigma_{P(\xi)=(\sigma,\tau),|\xi|=|\sigma|+|\tau|}[\xi]$ (the vertex of the chain complex $[\xi]$ is in order), we have\\
\hspace*{16mm}$(E^2_{*,*,*},d_2)=(\oplus_{P(\xi)=(\sigma,\tau),|\xi|=|\sigma|+|\tau|,\xi\neq\phi}
[\w H^*({\rm link}_K\sigma;\ff){\otimes}\w H^*({\rm link}_L\tau;\ff)]_{\xi},d_2)$.\\
So the correspondence $[x]_{\sigma}\otimes[y]_{\tau}\to \Sigma_{P(\xi)=(\sigma,\tau),|\xi|=|\sigma|+|\tau|}[x{\otimes} y]_{\xi}$ for all $x\in \w H^*({\rm link}_K\sigma;\ff)$ and $y\in \w H^*({\rm link}_L\tau;\ff)$
is an isomorphism from $(L_{*,*}(K;\ff),\partial_1){\otimes}(L_{*,*}(L;\ff),\partial_1)$ to $(E^2_{*,*,*},d_2)$. So $E^3_{*,*,*}=LH_{*,*}(K;\ff){\otimes}LH_{*,*}(L;\ff)$. But the spectral sequence collapse from $r\geqslant 3$. The theorem holds.

\end{document}